\documentclass[12pt]{article}
\usepackage{graphicx}
\usepackage{amsfonts,amssymb, amsmath, theorem}
\usepackage[latin1]{inputenc}
\usepackage{graphics}
\usepackage{cite}

\newtheorem{Th}{{\mbox{$\;\;\;\;\;\,$}}Theorem}[section]

\newtheorem{Lem}{{\mbox{$\;\;\;\;\;\,$}}Lemma}[section]

\newtheorem{Rem}{{\mbox{$\;\;\;\;\;\,$}}Remark}[section]

\newcommand{\qed}{\hbox to 0pt{}\hfill$\rlap{$\sqcap$}\sqcup$}

\newenvironment{proof}{{\noindent\it\underline{Proof}}:}{\hfill{$\blacksquare$}}

\newenvironment{proofES}{{\noindent\it\underline{Proof of Theorem \ref{super}}}:}{\hfill{$\blacksquare$}}

\newenvironment{proofMS}{{\noindent\it\underline{Proof of Theorem \ref{superlineal}, case 1}}:}{\hfill{$\blacksquare$}}

\begin{document}
\title{Existence and multiplicity of periodic solutions for a hematopoiesis model} 

\date{}

\author{\small{Pablo Amster and Roc\'io Balderrama\footnote{This work was fully supported by the projects UBACyT 20020120100029BA and PIP 11220090100637 CONICET}}}

\date{}
\maketitle

{\centerline{\small\textit{Departamento de Matem\'atica, Facultad de Ciencias Exactas y Naturales}}}

{\centerline{\small\textit{Universidad de Buenos Aires \& IMAS-CONICET - Argentina.}}}

\medskip


\begin{abstract}

A general nonautonomous Mackey-Glass equation for the regulation of the hematopoiesis with several non-constant delays is studied. 
Using topological degree methods, we prove the existence and multiplicity of positive periodic solutions. 

\end{abstract}

\section{Introduction}  
The following nonlinear autonomous delay differential equation was proposed
by Mackey and Glass \cite{MG}  to study the regulation of hematopoiesis: 
\begin{equation}
\label{MG} 
\frac{dP(t)}{dt}=\frac{\lambda \theta^{n}P(t-\tau)}{\theta^n + P^{n}(t-\tau)}-\gamma P(t).
\end{equation}
Here $\lambda, \theta, n, \gamma, \tau$ are  positive constants, $P(t)$ is the concentration
of cells in the circulating blood and the flux function $f(v)=\frac{\lambda \theta^{n}v }{\theta^n + v^{n}}$  
of cells into the blood stream depends on the cell concentration at an
earlier time. The delay $\tau$ describes the time  between the start of cellular production in the 
bone marrow and the release of mature cells into the blood. It is assumed that the cells are lost at a rate 
proportional to their concentration, namely $\gamma P(t)$, where $\gamma$ is the decay rate.
 This equation is a model of `dynamic disease'. 
This type of population dynamics equations attracted the interest of many 
researchers. Different aspects and properties of (\ref{MG}) 
have been studied by various authors, see for example
 \cite{GBL,GTB,M}.

Most often, the environment is not temporally constant; thus, 
it is intuitive to assume that this fact influences many biological 
dynamical systems and suggests the need of considering 
time-dependent parameters. 
Moreover, as remarked in \cite{CH,LK,SA}, more realistic models are those in which 
 periodicity of the environment and time delay play a role
(for more details, see e.g. \cite{N}).
In view of this, the following model is proposed in \cite{SA}:
\begin{equation}
\label{S-A}
x'(t)=\frac{q(t)x(t)}{r+x^{n}(t-mT)}-p(t)x(t)
\end{equation}
where $m,n$ are positive integers, $p, q$ 
are positive $T-$periodic functions
and the delay $\tau:= mT$ is a multiple of the period determined by the environment.


In order to establish a more realistic model, it is convenient to 
introduce a more general delay that extends the two above-referred cases.
Instead of assuming that the delay is constant or a multiple of the 
period of the environment, more general models are obtained by 
assuming that the time delay $\tau$ is an arbitrary continuous nonnegative
$T-$periodic function depending on $t$. 
The more general equation 
\begin{equation}
\label{MGcont}
x'(t)=\frac{a(t)x(t-\tau(t))}{1+x^{n}(t-\tau(t))}-b(t)x(t)
\end{equation}
where 
$a,b$ and $\tau$ are continuous positive $T-$periodic functions
was studied for example in
\cite{WJ, WJX,WL,WLZ,YFW}. 
Different aspects of equation (\ref{MGcont}) have been considered; 
in particular, existence of positive $T$-periodic solutions was 
proven, in most cases, using appropriate fixed point theorems.
In \cite{WLZ}, coincidence degree theory was
employed to prove the existence of a positive $T-$periodic solution 
under a condition that can be regarded as a particular application 
of Theorem \ref{sub}, case $(2)$ below. 
Moreover, when $a(t)=\gamma b(t)$ for some 
$\gamma>0$ and when $\tau, a$ and $b$ are constant, the conditions 
$\gamma>1$ and $a>b$ respectively are both necessary and sufficient for the existence of positive $T$-periodic solutions.

The following more general model was studied in 
\cite{BB} and \cite{LYZ}:

\begin{equation}
\label{BB}
x'(t)=\sum_{k=1}^{M}\frac{r_{k}(t)x^{\delta}(t-g_{k}(t))}{1+x^{\gamma}(t-g_{k}(t))}-b(t)x(t).
\end{equation}
Here, $\gamma$ is a positive constant 
and $r_k, b$ are positive $T-$periodic continuous functions.
For $\delta= 1$, 
existence and uniqueness of positive $T$-periodic solutions was 
studied in \cite{BB} for the particular case of constant 
proportional delays $g_k\equiv l_k T$; 
moreover, for general continuous, 
positive $T-$periodic $g_k$, attractiveness of some specific positive periodic
solutions was studied. 
For the case $\delta=0$ and $g_k$ continuous positive and 
$T-$periodic, existence and uniqueness of positive 
 $T-$periodic solutions of (\ref{BB}) was proven 
in \cite{LYZ} by fixed point methods, provided that 
one of the following conditions
is satisfied:
 
 $$(1)\,\, \gamma\le 1\,\,\, \hbox{ or }\,\,\, (2)\,\, \gamma >1 \hbox{ and } 
(\gamma-1)\left(\frac{e^{\int_{0}^{T}b(u)du}}{e^{\int_{0}^{T}b(u)du}-1}\int_{0}^{T}\sum_{k=1}^{M}r_k(t)dt\right)^{\gamma}\le 1.$$


Motivated by the preceding discussion, we consider the 
following more general nonlinear nonautonomous model with several delays
\begin{equation}
\label{BBg}
x'(t)=\sum_{k=1}^{M}
\lambda_{k} r_{k}(t)\frac{x^{m_{k}}(t-\tau_{k}(t))}{1+x^{n_{k}}(t-\mu_{k}(t))}-b(t)x(t)
\end{equation}
where $r_{k}(t), b(t), \tau_{k}(t)$ and $\mu_{k}(t)$ are positive and $T$-periodic functions and 
$\lambda_{k}, m_k, n_k$ are positive constants.

Existence of solutions of \eqref{BBg} under appropriate conditions follows 
from several abstract results, although multiplicity results are more scarce. 
For example, in \cite{HW} and \cite{YFW} Krasnoselskii type 
fixed point theorem in cones were employed in order to obtain conditions 
for the existence of at least two
$T-$periodic solutions of the general equation 
\begin{equation}
\label{YFW}
x'(t)=-a(t)x(t)+f(t,x(t-\tau_1(t)),\ldots,x(t-\tau_n(t))).
\end{equation}
It is observed, however, that these results can be applied 
only to few particular sub-cases of (\ref{BBg}). 
In such cases, the conclusions are comparable to our results below.
Moreover, the existence of three nonnegative periodic solutions of (\ref{YFW})
was studied by using Leggett-Williams fixed 
point theorem in \cite{BX,PS, PSD,PSP}.
However, the conditions obtained in \cite{PS}, as pointed out by the authors,
are very difficult to apply to (\ref{BBg}) with $M=1, m=1, \tau=\mu$. 
Thus, they established a complementary result with more straightforward conditions that 
can be applied to this model. 
Unfortunately, 
in \cite{PSD}, the authors observed that this latter result was incorrect. 

In section \cite[Applications]{PSD}, the hematopoiesis model (\ref{BBg}) for $M=1, \tau=\mu$ was studied, and the conditions obtained by the authors
are similar to the ones proposed in Theorem $4.2$ $(1)$ below. 
It is worth mentioning that, in the referred result of \cite{PSD}, only 
two of the three $T$-periodic solutions are positive and the 
third one is positive if $f(t,0)$ is not identically zero. 
This assumption is very restrictive and clearly is not fulfilled in 
(\ref{BBg}). 
Moreover, all the mentioned works
do not contemplate the {\sl superlinear} case
of (\ref{BBg}) (that is, $m_k>n_k +1 $ for some $k$). 
Finally, we may mention the results in \cite{ZZB}, 
in which the existence of at least $2n$ solutions of (\ref{YFW}) is proven,
although the conditions are not applicable to our model.

Our goal in this paper is to establish sufficient criteria to guarantee, 
on the one hand, the existence of positive $T$-periodic solutions of (\ref{BBg}) and, on the other hand, 
multiplicity of such solutions. 
Using degree theory, we shall obtain a set of natural and easy-to-verify conditions for the 
existence of one or more solutions.


The following notation will be used throughout the paper. Let 
$$C_T:= \left\{u(t)\in C(\mathbb{R},\mathbb{R}): u(t+T)=u(t) \hbox{ for all } t \right\} $$
denote the space of continuous $T-$periodic functions and define, for $r<s$,
$$X_{r}^{s}:=\left\{u(t)\in C_T:  r < u(t) < s \hbox{ for all } t\right\}.$$   
The closure of $X^{s}_{r}$ shall be denoted by $cl(X^{s}_{r})$. 
The average, the maximum value and the minimum value of an arbitrary function 
$\varphi\in C_T$ shall be denoted respectively by $\overline \varphi$, 
$\varphi_{max}$ and $\varphi_{min}$, namely

$$\overline \varphi:= \frac{1}{T}\int_{0}^{T} \varphi(t)\, dt, \qquad 
\varphi_{max}=\max_{[0,T]} \varphi(t),\qquad \varphi_{min}=\min_{[0,T]} \varphi(t).$$

In order to simplify some computations, we set $y(t):=\ln(x(t))$ and transform (\ref{BBg}) into the equivalent equation
\begin{equation}
\label{BBg1}
y'(t)=\sum_{k=1}^{M}\lambda_{k} r_{k}(t)\frac{e^{m_{k}y(t-\tau_{k}(t))-y(t)}}{1+e^{n_{k}y(t-\mu_{k}(t))}}-b(t).
\end{equation}

Finally we define, for convenience, the function $\phi :\mathbb{R}\to \mathbb{R}$ by
\begin{equation}
\label{BBg1phi}
\phi(\gamma):= 
\sum_{k=1}^{M}\lambda_{k} \overline r_{k} \frac{e^{(m_{k}-1)\gamma}}{1+e^{n_{k}\gamma}} - \overline b.
\end{equation}

The proof of our results shall be based on the continuation method. Specifically, we shall apply 
the following existence theorem, established in \cite{AI}.

\begin{Th}
\label{cont} Assume there exist constants $\gamma_1<\gamma_2$ such that
\begin{enumerate}
 \item If $y\in cl(X^{\gamma_2}_{\gamma_1})$ satisfies
 \begin{equation}
  \label{BBg1sigma}
 y'(t)= \sigma\left(\sum_{k=1}^{M}\lambda_{k} r_{k}(t)\frac{e^{m_{k}y(t-\tau_{k}(t))-y(t)}}{1+e^{n_{k}y(t-\mu_{k}(t))}}-b(t)\right)
\end{equation}
 for some $\sigma \in (0,1]$,  
 then $y \in X^{\gamma_2}_{\gamma_1}$. 
 \item $\phi(\gamma_1) \phi(\gamma_2)<0$. 
\end{enumerate}
Then (\ref{BBg1}) has at least one solution in $X_{\gamma_1}^{\gamma_2}$.
\end{Th}

Roughly speaking, if $\phi$ has different 
signs at both ends of some interval $[\gamma_1, \gamma_2]\subset\mathbb{R}$ then the continuation theorem
guarantees the existence of a $T$-periodic solution $y$ of (\ref{BBg1}) such that $y(t)\in 
(\gamma_1, \gamma_2)$ for all $t$. 
However, the first condition of 
Theorem \ref{cont} 
requires, in some sense, that the sign of $\phi$ does not change too fast. 

The main part of our analysis shall be based on a study of the behavior of $\phi$. 
For the existence of solutions it suffices, in most cases, 
to consider its behavior at $\pm \infty$; for the multiplicity results a more careful study 
is needed, in order to find ``large enough'' intervals of positivity and negativity of $\phi$. With that 
end in mind, we shall consider the sets 
$$
M_1:= \{k:0<m_k<1\}, \; 
M_2:= \{k:m_k=1\}, \;
M_3:= \{k:1<m_k<n_k+1\}
$$
$$M_4:= \{k:m_k=n_k+1\},\; M_5:=\{k:m_k>n_k+1\}$$
and 
$$\phi_{i}(\gamma):=\sum_{k\in M_i}\lambda_{k} \overline r_{k} \frac{e^{(m_{k}-1)\gamma}}{1+e^{n_{k}\gamma}},$$
so we may write $\phi(\gamma)=\sum_{i=1}^{5}\phi_{i}(\gamma)-\overline b.$ For notation convenience, we also define $C:=T\overline b=\int_0^Tb(t)\, dt$. 

This setting proves to be useful, since the limits $\lim_{\gamma\to\pm\infty}\phi_i(\gamma)$ are easy to compute and, 
moreover, $\phi_i(\gamma)$ is strictly monotone for $i\neq 3$ and a sum of one-hump functions for $i=3$. Thus, the behavior of $\phi$ can be 
understood by studying the interaction of these different terms.

The paper is organized as follows. In the next section, we adapt the abstract continuation theorem \cite[Thm 2.1]{AI} to equation (\ref{BBg}) 
and prove the 
existence of positive $T$-periodic solutions for the different cases. 
In Section \ref{mult}, we give sufficient conditions for the existence 
of $2$, $3$ or $4$ positive $T$-periodic solutions. Finally, in Section \ref{exa} we give an example with at least $6$ positive $T$-periodic solutions. 

\section{Existence of positive $T-$periodic solutions.}

In order to present our existence 
results in a more comprehensive way, 
we shall consider three different cases: 
the {\sl superlinear} case ($m_{k}>n_{k}+1$ for some $k$ ), the the {\sl sublinear} case ($m_k<n_{k}+1$ for all $k$)
and the {\sl asymptotically linear} case ($m_k\le n_k + 1$ for all $k$ and $m_j=n_{j}+1$ for some $j$). 
We give a detailed proof only of the first result, since the other two follow similarly.


\begin{Th}
\label{super}
Assume $m_{j}>n_{j}+1$ for some $j$. 
Furthermore, assume that one of the following conditions is fulfilled:
\begin{enumerate}
\item  $m_k > 1$ for all $k$.
\item  $m_k \geq 1$ for all $k$, $m_i=1 $ for some $i$ and 
$\sum_{k\in M_2}\lambda_{k}r_{k}(t)e^{C}
<b(t)$ for all $t$.
\item $m_i < 1$ for some $i$ and $\sum_{k=1}^{M}\lambda_k r_{k}(t)\frac{e^{(m_{k}-1)\gamma_1}e^{m_{k}C}}{1+e^{n_{k}\gamma_1}}< b(t)$ for all $t$ 
and some constant $\gamma_1$.
\end{enumerate}
Then (\ref{BBg}) admits at least one positive $T-$periodic solution.
\end{Th}



\begin{Th}
\label{sub}
Assume $m_{k}<n_{k}+1$ for all $k$. Furthermore, assume that one of the following conditions is fulfilled:

\begin{enumerate}
\item  $m_i<1$ for some $i$.

\item $m_k\ge 1$ for all $k$, $m_i=1$ for some $i$ 
 and $\sum_{k\in M_2} \lambda_k{r_k(t)} > {b(t)}$ for all $t$.
 
\item  $m_k > 1$ for all $k$ and 
$$\sum_{k=1}^M \lambda_k{r_k(t)} \frac{e^{(m_k-1)\gamma_1}}{1+e^{n_k(\gamma_1+C)}}  > {b(t)}$$ 
for all $t$ and some arbitrary constant $\gamma_1$.
\end{enumerate}
Then (\ref{BBg}) admits at least one positive $T-$periodic solution.
\end{Th}


\begin{Th}
\label{asint}
Assume $m_{k}\le n_{k}+1$ for all $k$ and $m_j=n_j +1$ for some $j$.
Furthermore, assume that one of the following conditions is fulfilled:

\begin{enumerate}
\item $m_k>1$ for all $k$ and $\sum_{k\in M_4}\lambda_k r_{k}(t)e^{-C m_k}>b(t)$
for all $t$.
\item  $m_k\geq 1$ for all $k$, $m_i=1$ for some $i$,
$\sum_{k\in M_4}\lambda_k r_{k}(t)e^{-C m_k}>b(t)$
and $\sum_{k\in M_2}\lambda_k r_{k}(t)e^C <b(t)$
for all $t$. 
\item $0 < m_i<1$ for some $i$ and $\sum_{m_k \in M_4}\lambda_k r_{k}(t)e^{C n_k}<b(t)$
for all $t$. 
\end{enumerate}
Then (\ref{BBg}) admits at least one positive $T-$periodic solution.
\end{Th}


\begin{proofES}
Let $y$ be a $T-$periodic solution of (\ref{BBg1sigma}) 
with $0<\sigma\le1$, then
$y'(t)\geq -b(t)$ and hence 
$y(t_1)-y(t_2)\le\int_{0}^{T}b(t) dt$ 
for any $t_1\le t_2\le t_1 + T$. This implies, since $y(t)$ is $T-$periodic, that
$y_{max}-y_{min}\le\int_0^T b(t)dt=C$. 
Moreover, since $m_{k}>n_{k}+1$ for some $k$ it follows that $\phi(\gamma)>0$ when $\gamma$ is large enough. 
Assume that $y_{max}$ is achieved at some value $t^*$, then 

$$b(t^*)e^{y_{max}}=\sum_{k=1}^{M}\lambda_{k}r_{k}(t^*)
\frac{e^{m_k y(t^{*}-\tau_{k}(t^*))}}{1+e^{n_{k}y(t^{*}-\mu_{k}(t^*))}}
$$
$$\geq \sum_{k=1}^{M}\lambda_{k}r_{k}(t^*)\frac{e^{m_k (y_{max}-C)}}{1+e^{n_k y_{max}}}$$
and consequently
$$b(t^*)\geq\sum_{k=1}^{M}\lambda_{k}r_{k}(t^*)\frac{e^{(m_{k}-1)y_{max}}e^{-C m_k}}{1+e^{n_{k}y_{max}}}.$$
Again, since $m_{k}>n_{k}+1$ for some $k$ we deduce that
$y_{max}$ cannot be too large. 
Thus, we may fix $\gamma_2\gg 0$ such 
that $y_{max}<\gamma_2$ for every $y\in C_T$ satisfying (\ref{BBg1sigma}) 
and $\phi(\gamma_2)>0$. 
In a similar fashion, we look for $\gamma_1<\gamma_2$ such that $\phi(\gamma_1)<0$ and $y_{min}\neq \gamma_1$. 

\begin{itemize}
\item[{\bf Case 1: }] \textit{$m_k > 1$ for all $k$.}
Here $$\phi(\gamma) \to -\overline{b}\qquad \hbox{as 
$\gamma\to-\infty$}. $$ 
Let $y\in C_T$ be a solution of (\ref{BBg1sigma}) and fix $t_*$ such that 
$y(t_*)=y_{min}$, then
$$b(t_*)\le \sum_{k=1}^{M}\lambda_{k}r_{k}(t_*)\frac{e^{(m_k -1)y_{min}}e^{C m_k}}{1+e^{n_k y_{min}}}$$
Suppose that $y_{min}=\gamma_1$, then
$$b(t_*)\le\sum_{k=1}^{M}\lambda_k r_k(t)\frac{e^{(m_k-1)\gamma_1+Cm_k}}{1+e^{n_k \gamma_1}}.$$
The right-hand side of the latter inequality tends to zero as $\gamma_1\to -\infty$.
We deduce that $y_{min}$ cannot take arbitrarily large negative values; hence, it suffices to take $\gamma_1\ll 0$.

\item[{\bf Case 2}.] \textit{$m_k \geq 1$ for all $k$ and $m_j=1 $ for some $j$.}
In this case, 
 $$\phi(\gamma) \to \sum_{k\in M_2}\lambda_{k}\overline{r_{k}}-\overline{b}<
 \sum_{k\in M_2}\lambda_{k}\overline{r_{k}}e^C-\overline{b}<0$$
 as 
$\gamma\to -\infty$. 
On the other hand, if $y\in C_T$ satisfies 
(\ref{BBg1sigma}) 
then
$$b(t_*)\geq \sum_{k\in M_2}\lambda_{k}r_{k}(t_*)\frac{e^{(m_k-1)y_{min}+m_k C}}{1+e^{n_k (y_{min}+C)}}=
\sum_{k\in M_2}\lambda_{k}r_{k}(t_*)\frac{e^{C}}{1+e^{n_k (y_{min}+C)}}
$$
and, again, we deduce that $y_{min}$ cannot take too large negative values. 
Thus, it suffices to take $\gamma_1\ll0$.

\item[{\bf Case 3}.] \textit{$m_k < 1$ for some $k$.}   
From the hypothesis,
$$\phi(\gamma_1)=\sum_{k=1}^{M}\lambda_{k}\overline{r_{k}}
 \frac{e^{(m_{k}-1)\gamma_1}}{1+e^{n_{k}\gamma_1}}-\overline{b}\le
 \sum_{k=1}^{M}\lambda_{k}\overline{r_{k}}
 \frac{e^{(m_{k}-1)\gamma_1}e^{m_{k}C}}{1+e^{n_{k}\gamma_1}}-\overline{b}<0.$$
Moreover, if $y_{min}$ is achieved at some value $t_*$, then
 $$b(t_*) \le \sum_{k=1}^{M}\lambda_k r_{k}(t_*)\frac{e^{(m_{k}-1)y_{min}}e^{m_k C}}{1+e^{n_{k}y_{min}}}.$$
We conclude that $y_{min}\neq \gamma_1$.
\medskip

\end{itemize}
\end{proofES}

\begin{Rem}

It is easy to verify that that the second condition in Theorem \ref{asint} can 
be replaced by  
\smallskip

\noindent 2'. $m_k\geq 1 \hbox{ for all   } k, m_i=1 \hbox{ for some } i, 
\sum_{k\in M_4}\lambda_k r_k (t) e^{C n_k  }<b(t)$
and $
\sum_{k\in M_2}\lambda_k r_k (t) >b(t) \hbox{ for all t.}$
\end{Rem}


\section{Multiplicity}

\label{mult} 
In this section, we shall employ Theorem \ref{cont} to prove the existence of multiple solutions. It is worth noticing that, when
$\phi$ is monotone, it changes sign at most once and the method cannot be applied. 
When $\phi$ is not monotone, it is not enough to obtain intervals of positivity and negativity: 
as mentioned, it is required that $\phi$ does not change sign too rapidly. 
For a more detailed analysis, the following functions shall 
be helpful:

$$\alpha(\gamma,t):=\sum_{k=1}^{M}\lambda_k r_{k}(t)
      \frac{e^{(m_{k}-1)\gamma}e^{-C m_{k}}}{1+e^{n_k (\gamma+C)}}-b(t) $$
      
%
$$\beta(\gamma,t):=\sum_{k=1}^{M}\lambda_k r_{k}(t)
\frac{e^{(m_{k}-1)\gamma}e^{C m_{k}}}{1+e^{n_k (\gamma-C)}}-b(t)$$

As before, our results shall be split in three different theorems, for the superlinear, 
sublinear and asymptotically linear cases. More concretely, 
 

\begin{Th}
\label{superlineal}
Assume that $m_j >n_j +1$ for some $j$.

\begin{enumerate}
\item Let $m_k >1$ for all $k$ and $1< m_i <n_{i}+1$ for some $i$. Assume there exist 
constants  $\gamma_1<\gamma_2$ such that 
$$\alpha(\gamma_1,t)>0 > \beta(\gamma_2,t) \hbox{ for all $t$.}$$
Then (\ref{BBg}) admits at least $3 $ positive $T-$periodic solutions.

\item Let $m_k\geq 1$ for all $k$, $m_i =1$ for some $i$, $m_k\notin (1,n_{k}+1)$ for all $k$.
Assume that 
$$\sum_{k\in M_2}\lambda_k r_{k}(t)> b(t) \hbox{ for all t }$$
and there exists $\gamma_1$ such that 
  $$\beta(\gamma_1,t)<0 \hbox{ for all $t$.}$$
Then (\ref{BBg}) admits at least $2$ positive $T-$periodic solutions.

\item Let $m_k\geq 1$ for all $k$, $m_i =1$ for some $i$ and $1< m_s< n_{s}+1$ for some $s$.
Assume that 
$$\sum_{k\in M_2}\lambda_k r_{k}(t)e^C< b(t)\hbox{ for all t}$$
and 
there exist constants $\gamma_1<\gamma_2$
such that 
$$\alpha(\gamma_1,t)>0 > \beta(\gamma_2,t) 
\hbox{ for all t. }$$
  	   
Then (\ref{BBg}) admits at least $3$ positive $T-$periodic solutions.

\item Let $m_{i}<1$ for some $i$, $m_k\notin (1,n_{k}+1)$ for all $k$.
Assume   
 
$$\beta(\gamma_1,t)<0 \hbox{ for all t and some constant $\gamma_1$. }$$

Then (\ref{BBg}) admits at least $2$ positive $T-$periodic solutions.

\item Let $m_{i}<1$ for some $i$, $1<m_s < n_{s}+1$ for some $s$. Assume there exist 
 some constants $\gamma_1<\gamma_2<\gamma_3$ such that

$$\alpha(\gamma_2,t)>0 >\beta(\gamma_i,t) \hbox{ for all t, }$$
for $i=1,3$. Then (\ref{BBg}) admits at least $4$ positive $T-$periodic solutions.  	   

\end{enumerate}

\end{Th}

 
\begin{Th}
\label{sublineal}
Assume that $m_k <n_k +1$ for all $k$. 
\begin{enumerate}
\item Let $m_{k}>1$ for all $k$ and assume there exists a constant $\gamma_1$ such that
$$\alpha(\gamma_1,t)>0 \hbox{ for all t.}$$
Then (\ref{BBg}) admits at least $2$ positive $T-$periodic solutions. 
       
\item Let $m_k\geq 1$ for all $k$, $m_{i}=1$, $m_j>1$ for some $i,j$. Assume that
$$\sum_{k\in M_2}\lambda_k r_k(t)e^C<b(t) \hbox{ for all t }$$ 
and there exists a constant $\gamma_1$ such that
$$\alpha(\gamma_1,t)>0 \hbox{ for all t}.$$
Then (\ref{BBg}) admits at least $2$ positive $T-$periodic solutions.

\item Let $0<m_i <1$, $m_j>1$ for some $i,j$. Assume there exist
some constants $\gamma_1<\gamma_2$  such that      
$$\alpha(\gamma_2,t)>0 >\beta(\gamma_1,t)\hbox{ for all t. }$$
Then
(\ref{BBg}) admits at least $3$ solutions.

\end{enumerate}
\end{Th}


\begin{Th}
\label{lineal}
Assume that $m_k \le n_k +1$ for all $k$ and
$m_j =n_j +1$ for some $j$. 

\begin{enumerate}
\item Let $m_k >1$ for all $k$ and $1< m_i<n_i+1$ for some $i$. Assume that 
$$\sum_{k\in M_4}\lambda_k r_k (t)e^{C n_k}<b(t) \hbox{ for all t }$$
and there exists a constant $\gamma_1$ such that 
$$ \alpha(\gamma_1,t)>0 
  \hbox{ for all t}.$$
Then (\ref{BBg}) admits at least $2$ positive $T-$periodic solutions.
       
\item Let $0<m_i< 1$ for some $i$, $m_k \notin (1,n_k +1)$ for all $k$. Assume that       
$$  \sum_{k\in M_4}\lambda_k r_k (t)e^{-C m_k}>b(t) \hbox{ for all t}$$
and there exists $\gamma_1$ such that 

$$\beta(\gamma_1,t)<0 \hbox{  for all t. } $$
Then (\ref{BBg}) admits at least $2$ positive  $T-$periodic solutions.
            
\item Let $0<m_i< 1$ and $1< m_s <n_s +1$ for some $i,s$. Assume that
$$\sum_{k\in M_4}\lambda_k r_k (t)e^{C n_k }<b(t) $$
and there exist constants $\gamma_1<\gamma_2$ such that
$$\alpha(\gamma_2,t)>0> \beta(\gamma_1,t)  \hbox{ for all t.}$$

 Then (\ref{BBg}) has at least $3$ positive $T-$periodic solutions.
            
\end{enumerate}

\end{Th}

As before, we shall only prove the first case of Theorem \ref{superlineal}, since all the remaining cases follow in an analogous way.

\medskip

\begin{proofMS}
We shall apply 
Theorem \ref{cont} on open bounded sets $X_{\gamma_0}^{\gamma_1}$,
 $X_{\gamma_1}^{\gamma_2}$ and $X_{\gamma_2}^{\gamma_3}$, with $\gamma_0<\gamma_1$ and $\gamma_3>\gamma_2$ to be determined.  
To begin, observe that 
$$ \phi(\gamma)\to -\overline{b} \qquad \hbox{ as } \gamma \to -\infty $$
 and 

$$\phi(\gamma)\to +\infty \qquad \hbox{ as } \gamma \to+\infty. $$

In the same way of Theorem \ref{super} it is proven that, if $\gamma_0\ll 0$ then there exists $y\in X_{\gamma_0}^{\gamma_1}$ solution of (\ref{BBg1}).

On the other hand, for all $t$ it is seen that
$$\phi(\gamma_1)>\alpha(\gamma_1,t)>0.$$ 
Moreover, if $y\in cl(X_{\gamma_1}^{\gamma_2})$ is a solution
of (\ref{BBg1sigma}) with $0<\sigma\le1$  and $y_{min}=y(t_*)$, then

$$b(t_*)e^{y_{min}}=\sum_{k=1}^{M}\lambda_{k}r_{k}(t_*)
\frac{e^{m_k y(t_{*}-\tau_{k}(t_*))}}{1+e^{n_{k}y(t_{*}-\mu_{k}(t_*))}}
$$

$$>\sum_{k=1}^{M}\lambda_{k}r_{k}(t_*)
\frac{e^{m_k y_{min}}}{1+e^{n_{k}(y_{min}+C)}}>
\sum_{k=1}^{M}\lambda_{k}r_{k}(t_*)
\frac{e^{m_k y_{min}}e^{-C m_k}}{1+e^{n_{k}(y_{min}+C)}}.
$$
It follows that $y_{min}\neq\gamma_1$.

Furthermore, $$\phi(\gamma_2)<\beta(\gamma_2,t)<0$$
for all $t$ and we deduce as before that $y_{max}\neq \gamma_2$.

Finally, the existence of 
$\gamma_3\gg 0$ such that the problem has a solution 
$y \in X_{\gamma_2}^{\gamma_3}$ follows as in Theorem \ref{super}.
\end{proofMS}

\bigskip
The following lemma shows, in the context of Theorem \ref{superlineal} (case 1), 
that if $r_k, m_k$ and $n_k$ are given, then it is possible to find parameters 
$\lambda_k$ such that assumptions are fulfilled.
Analogous arguments are valid for the remaining cases. 

\begin{Lem}
\label{lemaXXt}
Let $r_k, b: \mathbb{R}\to\mathbb{R}_{>0}$ be continuous and 
$T-$periodic functions and $m_k, n_k \in \mathbb{R}_{>0}$ such that 
$m_{k}>1$ for all $k$, $1<m_j<n_j +1$ for some $j$, $m_{i}>n_i +1$ for some $i$.
Then there exist $\lambda_{k}$ and 
$\gamma_1<\gamma_2$ such that 
$$\alpha(\gamma_1,t)>0 >\beta(\gamma_2,t) \hbox{ for all $t$}.$$
\end{Lem}

\medskip

 \begin{proof}
Using the sets $M_i$ as before, we may write $\alpha$ and $\beta$ as 
$$\alpha(\gamma,t)=\sum_{i=1}^{5}\alpha_i(\gamma,t)-b(t),\qquad
\beta(\gamma,t)=\sum_{i=1}^{5}\beta_i(\gamma,t)-b(t).$$
Observe that, for each $t\in [0,T]$ and $i=1,\ldots,5$, 
$\alpha_i(\cdot,t)$ and $\beta_i(\cdot,t)$ have the same qualitative behavior as $\phi_i$.

We begin by setting the parameters $\lambda_k\in M_3$.
For arbitrary $\gamma_1$, take $\lambda_k \in M_3$ large enough such that
$$\alpha_3(\gamma_1,t)=\sum_{k\in M_3}
\lambda_{k} r_{k}(t)\frac{e^{(m_{k}-1)\gamma_1}e^{-C m_k}}{1+e^{n_{k}(\gamma_1 +C)}}
-b(t)>0.
$$
For $\epsilon\in (0, b_{min})$, there exists $R>\gamma_1$ such that 
$$\beta_3(\gamma,t)=\sum_{k\in M_3}
\lambda_{k} r_{k}(t)\frac{e^{(m_{k}-1)\gamma}e^{C m_k}}{1+e^{n_{k}(\gamma -C)}}
<\sum_{k\in M_3}
\lambda_{k} r_{k}^{max}\frac{e^{(m_{k}-1)\gamma}e^{C m_k}}{1+e^{n_{k}(\gamma -C)}}<\epsilon$$
for $\gamma>R$ and all $t$. Thus, we may fix $\gamma_2 >R$ and proceed with the remaining parameters. 

Next, for $k\in M_4\cup M_5$ we set $\lambda_k$ small enough so that 
 $$\sum_{k\in M_4\cup M_5}
\lambda_{k} (r_{k})_{max}\frac{e^{(m_{k}-1)\gamma_2}e^{C m_k}}{1+e^{n_{k}(\gamma_2-C)}}
<b_{min}-2\epsilon$$
and hence
$$(\beta_4+\beta_5)(\gamma_2,t)=\sum_{k\in M_4\cup M_5}
\lambda_{k} r_{k}(t)\frac{e^{(m_{k}-1)\gamma_2}e^{C m_k}}{1+e^{n_{k}(\gamma_2-C)}}
<b_{min}-2\epsilon<b(t)-2\epsilon.
$$
Thus the conclusion follows since
$$\beta(\gamma_2, t)=
(\beta_3+\beta_4+\beta_5)(\gamma_2,t)-b(t)
<\epsilon-2\epsilon<0$$ 
and 
$$\alpha(\gamma_1,t)
=(\alpha_3+\alpha_4+\alpha_5)(\gamma_1,t)-b(t)$$

$$>\alpha_3(\gamma_1,t)-b(t)>0$$

\end{proof}

\section{Example}

\label{exa} 

As shown in Theorem \ref{superlineal}, case 5,
equation (\ref{BBg}) has at least $4$ positive $T$-periodic solutions.
The following example shows that, in fact, the problem may have more solutions.
Let $k=4$ and 
$b(t)=1.1+0.02\cos(\frac{2\pi t}{T})$, $T=0.005$,
 $m_1=0.95$, $n_1=2$, $\lambda_1r_1(t)=0.04+0.002\cos(\frac{2\pi t}{T})$,
 $m_2=4.73$, $n_2=3.74$, $\lambda_2r_2(t)=1.3+0.002\cos(\frac{2\pi t}{T})$,
 $m_3=1.0001$, $n_3=10.2$, $\lambda_3r_3(t)=0.9+0.002\cos(\frac{2\pi t}{T})$,
 $m_4=1.12$, $n_4=0.11$, $\lambda_4r_4(t)=0.06+0.002\cos(\frac{2\pi t}{T})$.

Set $\gamma_1=-5$, $\gamma_2=-0.3$ $\gamma_3=0.2$, $\gamma_4=5$, $\gamma_5=34$. It is verified (see Figure $1$) that
$$\alpha(\gamma_2,t)>0.09,\alpha(\gamma_4,t)>0.1 \hbox{ for all $t$ }$$
and
$$\beta(\gamma_1,t)<-0.08,\beta(\gamma_3,t)<-0.01,\beta(\gamma_5,t)<-0.01 \hbox{ for all $t$. } $$
Moreover, since  $0<m_1=0.95<1$ and 
$m_4=1.12>n_4+1=1.11$, it follows that 
$$\lim_{\gamma\to-\infty}\phi(\gamma)=
\lim_{\gamma\to+\infty}\phi(\gamma)=
+\infty.$$
Thus, we conclude that (\ref{BBg}) has at least six positive $0.005-$periodic solutions for arbitrary nonnegative $0.005$-periodic delays $\tau_k, \mu_k$. 

\begin{figure}[ht!]
\centering
\includegraphics[scale=0.45]{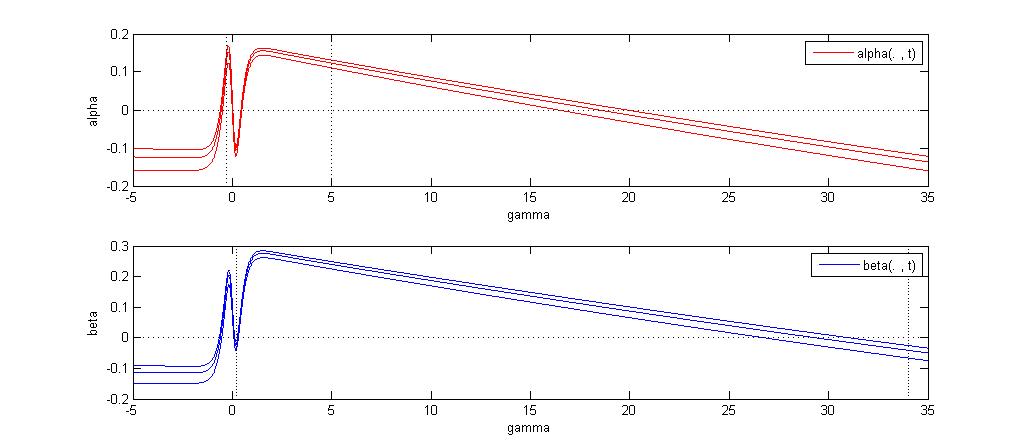}
\parbox{0.8\textwidth}{
\caption{\footnotesize $\alpha(\gamma,t)$ and $\beta(\gamma,t)$ for each $t\in[0:0.01:T].$  }}
\end{figure}

\end{document}